\documentclass[12pt]{article}
\usepackage{amsmath,amscd,amssymb}
\author{Theodore Voronov\thanks{{\it Current address:\/} Department of Mathematics, University of California, Berkeley, Berkeley, CA 94720. {\it Email:\/} {\tt theodore@math.berkeley.edu}}}
\title{Quantum categories. Quantization of the category of linear spaces.}
\date{}
\DeclareMathSymbol\leqslant{\mathrel}{AMSa}{"36} 
\DeclareMathSymbol\geqslant{\mathrel}{AMSa}{"3E} 
\numberwithin{equation}{section}
\begin{document}
\maketitle
\vspace{-1cm}
\begin{center}
\small       Moscow State University\\
       Department of Mechanics and Mathematics\\
      Chair of Higher Geometry and Topology\\
        Vorobyovy Gory,  Moscow 119899, Russia \\
       {\it Email:}\/  {\tt theodore@mech.math.msu.su}
\end{center}

%
%
\renewcommand{\a}{\alpha    }
\renewcommand{\b}{\beta    }
\renewcommand{\c   }{\gamma    }
\newcommand{\C   }{\Gamma    }
\renewcommand{\d   }{\delta    }
\newcommand{\D   }{\Delta    }
\newcommand{\e   }{\varepsilon}
\newcommand{\g   }{\gamma    }
\newcommand{\G   }{\Gamma    }
\newcommand{\h   }{\eta    }
\renewcommand{\l   }{\lambda    }
\renewcommand{\L   }{\Lambda    }
\newcommand{\m   }{\mu    }
\newcommand{\n   }{\nu    }
\renewcommand{\O   }{\Omega }
\renewcommand{\o   }{\omega    }
\newcommand{\p   }{\pi    }
\renewcommand{\P   }{\Pi }
\newcommand{\s   }{\sigma    }
\newcommand{\Ss   }{\Sigma    }
\renewcommand{\t   }{\tau    }
\renewcommand{\u   }{\upsilon    }
\newcommand{\U   }{\Upsilon    }
\renewcommand{\v   }{\upsilon    }
\newcommand{\V   }{\Upsilon    }
\newcommand{\x   }{\xi    }
\newcommand{\X   }{\Xi    }
\newcommand{\y   }{\upsilon    }
\newcommand{\Y   }{\Upsilon    }
\newcommand{\z   }{\zeta    }
\newcommand{\Z   }{\Zeta    }
\newcommand{\ps   }{\psi    }
\newcommand{\ph   }{\varphi    }
\newcommand{\PH   }{\varPhi    }
%

\begin{abstract}
We define generalized bialgebras and Hopf algebras
and on this basis we introduce quantum categories and quantum groupoids. The quantization of the
category of linear (super)spaces is constructed. We
establish a criterion for the classical value of the dimension of
the polynomial function algebra   on the full  subcategory of this quantum category 
specified by  Sudbery type commutation relations for quantum vector
spaces  (the ``Poincare -- Birkhoff -- Witt property'').
The criterion is the equality of the ``quantum constants''
$c_\a=c_\b^{\pm 1}$ of  the  quantum vector spaces.  (These constants appeared in earlier works on PBW for the multiparameter quantization of the general linear group.)
Links with the Yang--Baxter equation and the
 Yang -- Baxter structures of quantum linear spaces are established.
The role of categories as a generalization of groups is discussed.
\end{abstract}


\newcommand{\qed}{$\quad\square$}
\newcommand{\lra}{\longrightarrow}
\newcommand{\lla}{\longleftarrow}
\newcommand{\llra}{\longleftrightarrow}
\newcommand{\iso}{\stackrel{\sim}{=}}
\newcommand{\Hom}{\operatorname{Hom}}
\newcommand{\End}{\operatorname{End}}
\newcommand{\GLn}{\operatorname{GL}\,(n)}
\newcommand{\GLnm}{\operatorname{GL}\,(n\mid m)}
\newcommand{\Ann}{\operatorname{Ann}}
\newcommand{\Ztwo}{\mathbb Z_2}
\newcommand{\Linq}{{\mathsf L\mathsf i\mathsf n}_q}
\newcommand{\Linqs}{{\mathsf L\mathsf i\mathsf n}_q^{(s)}}
\newcommand{\Linqplus}{{\mathsf L\mathsf i\mathsf n}_q^+(\l)}

\newcommand{\Tab}{T_{\a \b}}
\newcommand{\Tac}{T_{\a \c}}
\newcommand{\Tbc}{T_{\b \c}}
\newcommand{\detab}{\operatorname{det}_{\a \b}}
\newcommand{\detac}{\operatorname{det}_{\a \c}}
\newcommand{\detbc}{\operatorname{det}_{\b \c}}

\newcommand{\Aa}{A_{\a}}
\newcommand{\Ab}{A_{\b}}
\newcommand{\Ac}{A_{\c}}
\newcommand{\Ba}{B_{\a}}
\newcommand{\Bb}{B_{\b}}
\newcommand{\Aaa}{A_{\a \a}}
\newcommand{\Aab}{A_{\a \b}}
\newcommand{\Aba}{A_{\b \a}}
\newcommand{\Abc}{A_{\b \c}}
\newcommand{\Aac}{A_{\a \c}}
\newcommand{\Aad}{A_{\a \d}}
\newcommand{\Abb}{A_{\b \b}}
\newcommand{\Abd}{A_{\b \d}}
\newcommand{\Acd}{A_{\c \d}}
\newcommand{\Mab}{M_{\a \b}}
\newcommand{\Maa}{M_{\a \a}}
\newcommand{\Mac}{M_{\a \c}}
\newcommand{\Mbc}{M_{\b \c}}

\newcommand{\Dabc}{\D_{\a \b \c}}
\newcommand{\Dbcd}{\D_{\b \c \d}}
\newcommand{\Dacd}{\D_{\a \c \d}}
\newcommand{\Daab}{\D_{\a\a\b}}
\newcommand{\Daba}{\D_{\a\b\a}}
\newcommand{\Dbab}{\D_{\b\a\b}}
\newcommand{\Dabb}{\D_{\a\b\b}}
\newcommand{\Dabd}{\D_{\a\b\d}}
\newcommand{\dab}{\d_{\a \b}}
\newcommand{\dbc}{\d_{\b \c}}
\newcommand{\dac}{\d_{\a \c}}
\newcommand{\ea}{\e_\a}
\newcommand{\eb}{\e_\b}

\newcommand{\qa}{q_\a}
\newcommand{\qb}{q_\b}
\newcommand{\qc}{q_\c}
\newcommand{\pa}{p_{\a}}
\newcommand{\pb}{p_{\b}}
\newcommand{\qab}{q^{AB}}
\newcommand{\qba}{q^{BA}}
\newcommand{\qaa}{q^{AA}}
\newcommand{\qkl}{q^{KL}}
\newcommand{\qkn}{q^{KN}}
\newcommand{\qlk}{q^{LK}}
\newcommand{\qln}{q^{LN}}
\newcommand{\qnl}{q^{NL}}
\newcommand{\pab}{p^{AB}}
\newcommand{\pba}{p^{BA}}
\newcommand{\paa}{p^{AA}}
\newcommand{\pkl}{p^{KL}}
\newcommand{\pkn}{p^{KN}}
\newcommand{\plk}{p^{LK}}
\newcommand{\pln}{p^{LN}}
\newcommand{\pnl}{p^{NL}}
\newcommand{\qabn}{q_{AB}}
\newcommand{\qacn}{q_{AC}}
\newcommand{\qbcn}{q_{BC}}
\newcommand{\qban}{q_{BA}}
\newcommand{\qcan}{q_{CA}}
\newcommand{\qbdn}{q_{BD}}
\newcommand{\qkln}{q_{KL}}
\newcommand{\pabn}{p_{AB}}
\newcommand{\pacn}{p_{AC}}
\newcommand{\pbcn}{p_{BC}}
\newcommand{\pban}{p_{BA}}
\newcommand{\pcan}{p_{CA}}
\newcommand{\pbdn}{p_{BD}}
\newcommand{\pkln}{p_{KL}}
\newcommand{\pv}{P_V}
\newcommand{\qv}{Q_V}
\newcommand{\pw}{P_W}
\newcommand{\qw}{Q_W}

\newcommand{\VIJ}{(V,I,J)}
\newcommand{\Vq}{V_q}
\newcommand{\VqP}{V_q\P}
\newcommand{\VI}{(V',I)}
\newcommand{\VPJ}{(\P V',J^\P)}

\newcommand{\At}{\tilde A}
\newcommand{\Bt}{\tilde B}
\newcommand{\Ct}{\tilde C}
\newcommand{\Dt}{\tilde D}
\newcommand{\Kt}{\tilde K}
\newcommand{\Lt}{\tilde L}

\newcommand{\tak}{t_A{}^K}
\newcommand{\tas}{t_A{}^S}
\newcommand{\sak}{s_A{}^K}
\newcommand{\tbl}{t_B{}^L}
\newcommand{\tbk}{t_B{}^K}
\newcommand{\tks}{t_K{}^S}
\newcommand{\tal}{t_A{}^L}
\newcommand{\tab}{t_A{}^B}
\newcommand{\tck}{t_C{}^K}
\newcommand{\tdl}{t_D{}^L}
\renewcommand{\tan}{t_A{}^N}
\newcommand{\tbm}{t_B{}^M}


\section*{Introduction}
\subsection*{ Motivation.}

The heuristic meaning  of quantum group theory consists in the consideration
of transformations
that depend on noncommuting parameters.
The parameters serve as generators of a certain (noncommutative)
algebra, which is interpreted as the
``function algebra" of the quantum group. The group nature of the
transformations is formalized in the assumption that this algebra is endowed
with the Hopf structure.
In matrix case the Hopf structure is induced just by the standard matrix
multiplication and by taking matrix inverse.
But the transformations of the group nature (and the related algebraic
structures) are not the only possible. Quantum semigroups appear along with
quantum groups (at least as the intermediate product, see
~\cite{rtf},~\cite{dem2},~\cite{dem3}). Multivalued Hopf algebras were defined
by Buchstaber and Rees~\cite{buch2}. For Novikov's operator doubles (which
are not Hopf algebras) the ``Hopf type'' questions can be studied
~\cite{nov},~\cite{buch1},~\cite{buch3}. Non-standard diagonals appear in the
Odessky -- Feigin algebras~\cite{of}. All this formally exceeds the
boundaries of quantum groups but nevertheless is in coherence with the
heuristic thesis stated above.

In the present paper we introduce and study  QUANTUM CATEGORIES and the
simplest but fundamental example, the QUANTIZATION OF THE CATEGORY OF
LINEAR (SUPER)SPACES.

This needs to be somehow explained since it is very common to treat
categories only as ``general nonsense''.
Quantum group theorists  use
categories mostly in this manner. (Categories of
representations and their abstract generalizations like ``tensor'' or
``braided'' categories are good example.) But the philosophy of categories
in this paper is
quite different. We consider categories first of all as algebraic
structures, sets of arrows with the multiplication law.
Hence categories are an
immediate generalization of groups, for many reasons more useful
than semigroups. Actually, this structure appears quite often in  very
classical situations. The category structure of such examples is
of course known but usually is not stressed.
(To mention just few.
Paths on a topological space form a category.  Parallel transports
in fiber bundles and the homology ``local coefficients'' are  examples of
its representations. Cobordism provide examples of categories with
``films'' as arrows. This includes  ``histories'' or Lorentzian
cobordism, see~\cite{mis},~\cite{gib}. It is worth mentioning that the
coordinate changes in the fixed open domain $U\subset\mathbb R^n$ also form a
category, not a group. Its representations exactly correspond to
``geometrical objects''.) Recent  research supplied new examples
of categories naturally appearing in the context very far from
``general nonsense''. For example, it was found that the representations of
classical groups, such as spinor representation, naturally extend to certain
categories of linear
relations~\cite{ner1}, ~\cite{ner2}. The study of  Poisson geometry lead to
symplectic groupoids~\cite{kar}, which provide an example of ``continuous
categories''. We can mention the category of tangles~\cite{tur}, etc.
This list can be extended.

In the context of ``general nonsense'', the role of  categories should be compared
to that of groups
in Klein's classification of geometries.
(This view was expressed by the founders of category theory themselves~\cite{em}.)
Hence the main thesis of a {\sl ``category program'' }\/
can be stated as follows.
{\it Categories are the most important generalization of groups; it is
necessary to study their actions, representations, deformations
and extensions,
cohomology} (see~\cite{bau},~\cite{wat}, ~\cite{qui}), {\it the analogues of
Lie theory}\/ etc.
In fact, 
the research in these directions is going on, with various
motivations.
We can especially note ~\cite{ner2} and~\cite{kar}.
(The author himself was brought to the necessity of a program
like this through the reflections on some problems of infinite-dimensional
geometry connected with the topics of~\cite{v2}.)
The study of  quantum categories initiated in the present paper
can be considered as one more step in the
realization of the ``category program''.

Inside
the quantum group theory the particular motivations for introducing
quantum categories are the following.
The same classical object can admit different ``quantum deformations'' (for
example, corresponding to different values of the deformation parameter).
Which
should be chosen ? The correct answer is that all of them should be
considered simultaneously, together with all possible transformations
between them. This inevitably leads to a quantum category or a quantum groupoid (if restricted to invertible transformations), 
while quantum groups  correspond to  particular
choices of deformation.

Let us consider a simple example. It is new and can serve as
a very good illustration of the
basic ideas of this work.

\medskip
{\bf Example.} Take matrix elements of
    $T_{\a \b }=\left(\begin{array}{cc} a&b \\
                                        c&d
                       \end{array}\right)$
as generators of an associative algebra $M_{\a \b}$ and impose
the following commutation relations
\begin{equation}
\begin{aligned}
    ab- \l \qb ^{-1}ba&=0, \\
    ac- \l \qa ca&=0,\\
    ad- \qa \qb ^{-1}da&=(\l - \l ^{-1})\qa cb\\
    bc -\qa \qb cb &=0,\\
    bd -\l \qa db&=0,\\
    cd - \l \qb ^{-1}dc&=0,
\end{aligned}
\end{equation}
where $\l \neq 0, \pm i$, \ $\qa, \qb \neq 0$.
(We shall treat $T_{\a \b}$
as the matrix of a  homomorphism from the quantum deformation of the linear
space $\mathbb C^2$ with the parameter $\qa$ to that with the parameter $\qb$.
This is explained in the main text below.)
The category nature of the transformation $T_{\a \b}$ is reflected in the
fact that if we take the matrix
$T_{\a \b}$ suiting (1) and the matrix $T_{\b \c}$ suiting
(1) with $\qb$ instead of $\qa$ and $\qc$ instead $\qb$,
then the product
 $T_{\a
\c}=T_{\a \b}T_{\b \c}$ also suits (1) with $\qa$, $q_\c$ as the
parameters
(the matrix elements of the different matrices are assumed to be
commuting).
If
$\qa=\qb=\l=1$ then we obtain the usual algebra of functions on the matrix
semigroup
$\mathrm{M}\mathrm{a}\mathrm{t}(2)$.
Thus it is wise to think that after the quantization this semigroup
becomes a ``quantum category"
(with $\qa$ as ``objects") and the group $\mathrm{G}\mathrm{L}(2)$
after the quantization
turns a ``quantum groupoid''
if one adds a formal inverse
$(\detab(\Tab))^{-1}$ to each of the algebras $\Mab$, where $\detab(\Tab)$
stands for $ad-\l \qa cb=ab-\l \qb^{-1}bc$.
The algebras $\Mab$ possess the ``Poincare - Birkhoff - Witt'' property:
the ordered monomials in matrix entries form the additive bases.
For $\qa=\qb=1$ the algebra $\Maa$ reduces to
the function algebra on
the quantum group $\mathrm{G}\mathrm{L}(2)_q$  defined in~\cite{rtf},
with $q=\l$.

\bigskip
(The appearance of ``quantum categories" in such  context can be related
to the thesis QUANTIZATION ELIMINATES DEGENERACY suggested by
Reshetikhin, Takhtajan, Faddeev \cite{rtf}.)

The results  of this paper were obtained in November 1994 -- January 1995.
They were reported at the seminar on quantum groups under guidance of
J.~Bernstein, J.~Donin and S.~Shnider in Bar-Ilan University (Ramat-Gan,
Israel).
\bigskip
\subsection*{ Contents.} The paper is divided into three sections.

In  Section 1 we
introduce  generalized bialgebras which are in the same relation to
quantum categories as  usual Hopf algebras are to quantum groups.

In  Section 2 the quantum deformation of the category of vector superspaces
is constructed. We make use of the method of the ``universal coacting'' known
for
quantum groups \cite{man1},\cite{man2},\cite{dem2},\cite{dem3}. We
calculate the commutational relations and consider various examples
including the quantization of the dual space, the quantization of
bilinear forms etc.

In  Section 3 we present a generalized
 \mbox{$R$-matrix} form for the commutational relations on
the quantum category
 $\mathsf{L}\mathsf{i}\mathsf{n}_q$ and establish a criterion for the
classical value of the dimension of the function algebra.
This is a technically hard theorem, which may be considered as the main
theorem of this paper. The connection with the ``Yang -- Baxter
structures'' on quantum linear spaces is established.

In the Appendix we give a general definition
of  ``indexed bialgebras'' which include
both the bialgebras defined in  Section 1 and those dual to them.

\subsection*{ Notation.}
Throughout the paper the standard tensor notation is used with the regard of
\mbox{$\mathbb Z_2$-grading}. (For the ``super'' notions that are used we
refer to
\cite{v1},\cite{v3}.) For some reason it is convenient to write the
coordinates of a vector as a row and those of a covector as a column,
and to consider their pairing in the form
$\langle v, v' \rangle,\, v \in V, \, v' \in V' $.

\subsection*{ Acknowledgements.}

I wish to thank V.M.~Buchstaber for the discussions of the results of this
paper and for supplying the texts of
~\cite{buch1},~\cite{buch2},~\cite{buch3}, E.E.~Demidov for supplying
the texts of~\cite{dem2} and~\cite{man1}, and O.M.~Khudaverdian for
the discussions of the role of categories in connection to mathematical
physics. I also wish to thank J.~Bernstein,
J.~Donin and S.~Shnider for the discussions in  Tel-Aviv and
Bar-Ilan universities and for their hospitality in May--June 1995. I thank
J.~Donin and S.~Shnider for the text of their paper~\cite{don} and I thank
S.~Shnider for the copy of his remarkable book with S.~Sternberg~\cite{ss}.

\renewcommand{\hom}{\operatorname{\mathbf{H}\mathbf{o}\mathbf{m}}}
\newcommand{\ssmooth}{\operatorname{\mathbf{C}^\infty}}
\newcommand{\Mor}{\operatorname{Mor}}
\newcommand{\mor}{\operatorname{\mathbf{M}\mathbf{o}\mathbf{r}}}
\newcommand{\Ob}{\operatorname{Ob}}
\newcommand{\Cat}{\mathsf{C}}

\section{Generalized bialgebras}

We are going to introduce a generalization of bialgebras as to make
possible the discussion of  ``quantized categories", in the same way as
the usual bialgebras and Hopf algebras serve as the foundation of the
notion of a ``quantum group". (More general definition of an ``indexed
bialgebra" will be discussed in the appendix.)

Consider a set of indices $\L$ (``objects'').

{\bf Definition 1.1.} A  set $A=(\Aab),\quad \a,\b\in\L$ of
linear spaces or modules over a commutative
ring
is said to be a {\it  coalgebra with objects
 }$\L$ if  linear maps $\Dabc : \Aac \lra
\Aab\otimes\Abc$ are given, for all $\a,\b,\c\in\L$.
It is called a
{\it coassociative coalgebra} if the diagram
\begin{equation}
\begin{CD}
          \Aad                @>{\Dabd}>>           \Aab \otimes \Abd \\
        @V{\Dacd}VV                               @VV{1 \otimes \Dbcd}V \\
     \Aac\otimes \Acd  @>>{\Dabc \otimes 1}>  \Aab\otimes \Abc\otimes \Acd
\end{CD}
\end{equation}
is commutative for every $\a,\b,\c,\d\in\L$.

A set of linear maps $\ea:\Aaa\lra R$
(here $R$ is the main field or ring) with the property that all
diagrams of the form
\begin{equation}
\begin{CD}
  \Aaa\otimes\Aab    @<{\Daab}<<    \Aab     @>{\Dabb}>>   \Aab\otimes\Abb \\
   @V{\ea\otimes1}VV                @VVV                 @VV{1\otimes\eb}V \\
    R\otimes\Aab          @>>>      \Aab       @<<<        \Aab\otimes R
\end{CD}
\end{equation}
must be commutative will be called a {\it  counit} in  $A$.
(Here  $R\otimes\Aab\lra\Aab, \Aab\otimes R\lra\Aab$
are the natural isomorphisms and middle vertical arrow is the identity
map.)

\medskip
The definition of the {\it homomorphism } of coalgebras
$f:(\Aab)_{\a,\b\in\L}\lra(B_{i,j})_{i,j\in I}$ is obvious.
(A map of sets $\a\mapsto i(\a)$ must be given.)

\medskip
{\bf Definition 1.2.} Let all modules be  the associative
$R$-algebras with units and all maps $\Dabc$ and $\ea$ be the
$R$-linear algebra homomorphisms. Then we shall call $A$ a {\it
bialgebra with objects} $\L$.

In the sequel we shall simply say ``coalgebra'' and ``bialgebra'' having in
mind the given definitions.

\medskip
{\bf Example 1.1.} {\sl Functions on a category.} Let $\Cat$ be a small
category with the object set
$\Ob{\Cat}=\L$. Define  $\Aab$ as
the function algebra on  $\Mor(\b,\a)$. Then the composition induces
the comultiplication $\Dabc:\Aac \lra \Aab\otimes\Abc$ and the embedding of
the unit $1_\a\in
\Mor(\a,\a)$ induces the counit homomorphism $\ea:\Aaa \lra R$.
(We assume the category and the class of functions in consideration to be
such that the Cartesian product of the morphism sets corresponds to a
tensor product of function algebras.)
We obtain the bialgebra with the object set
 $\L$. Here all algebras $\Aab$ are commutative.

\medskip
{\bf Example 1.2.} {\sl Additive categories.} An additive category
$\mathsf{A}$ can be treated as ``ring with several objects''
\cite{mit}. The additive category axioms are exactly dual
(in the sense of reverting arrows)
to our definition of a coalgebra  (over $\mathbb Z$,
coassociative and with a counit). An important case is
additive categories like $\mathbb Z \Cat$ or $R\Cat$, generated by the
morphisms
of an arbitrary small category $\Cat$ (an analogue of the group ring).
In this case
each module $\Aab=R\Mor_{\Cat}(\b,\a)$ is endowed with the comultiplication
$\Aab \lra \Aab\otimes\Aab$, defined on $a_i \in
\Mor_{\Cat}(\b,\a)$ by the formula $\D a_i=a_i\otimes\a_i$
(the same as for a  group
ring). Having this additional structure in mind, we may say that the
category algebra is a ``bialgebra'', but not in the sense of the
definition above but of the ``dual'' definition (w.r.t. the reverse arrows).
The comultiplication is symmetric here.
Considering the dual modules
$\Aab'=\Hom(\Aab,R)$ we (under certain assumptions) come to a bialgebra.
It coincides with the bialgebra of the previous example if, for example, all
sets $\Mor(\a,\b)$ are finite.

\medskip
{\bf Example 1.3.} ``{\it Supercategories}''. They can be considered in the
same fashion as Lie supergroups. (For instance one can obviously
define the supercategory
of $\Ztwo$-graded vector spaces, which contains the general linear
supergroups,
and in the same way the (infinite-dimensional) supercategory of smooth
manifolds.)
If for a given supercategory we consider function algebras on the
supermanifolds
$\mor(\b,\a)$, then, with standard comments on  tensor products, we again
obtain a bialgebra with the multiplication commutative in
$\Ztwo$-graded sense.

\medskip
Thus categories and supercategories can be described by bialgebras with
commutative multiplication. It is natural to think that arbitrary
bialgebras (in our sense) without commutativity condition may be associated
with  ``quantized categories''. Sure, it makes sense in the case  a
classical (super)category is at hand and the ``quantum category'' in
consideration reduces to it for certain values of the ``parameters''
(in most broad sense) on which the construction depends.

\medskip
{\bf Definition 1.3.} By an {\it antipode} in a bialgebra $A=(\Aab)$ we
shall mean a set of algebra antihomomorphisms
 $S_{\a\b}:\Aab \lra \Aba$ such that
the diagram
{
\unitlength=1em
\newcommand{\larr}{\begin{picture}(0,0)
                     \put(-3,0.28){\vector(1,0){5.8}}
                     \put(-1,0.7){${\scriptstyle 1\otimes S_{\b\a}}$}
                   \end{picture}}
\newcommand{\rarr}{\begin{picture}(0,0)
                     \put(2.8,0.28){\vector(-1,0){5.8}}
                     \put(-1,0.7){${\scriptstyle S_{\b\a}\otimes 1}$}
                   \end{picture}}
\begin{equation}
\begin{CD}
\Aab\otimes\Aba @. \larr   @.  \Aab\otimes\Aab  @. \rarr @. \Aba\otimes\Aab \\
 @A{\Daba}AA            @.    @VmVV               @.          @AA{\Dbab}A \\
 \Aaa         @>>{\ea}> R @>>> \Aab          @<<< R  @<<{\eb}< \Abb
\end{CD}
\end{equation}
is commutative. (The unmarked arrows are inclusion of unit.)
}

Bialgebras with  antipode correspond to
``quantum groupoids'', i.e. to quantum categories
in which ``all arrows are invertible''.

\medskip
In classical situation the notion of a ``concrete category''
is important (that is a subcategory of the category of sets).
A bialgebra  $(\Aab)$ is called  {\it
concrete}, if the ``coaction'' on a given set of algebras $(\Aa)$ is
provided.
That is the algebra homomorphisms $\dab:\Ab \lra \Aa\otimes\Aab$
that are compatible with
the comultiplication: $(\dab\otimes1)\circ\dbc=(1\otimes\Dabc)\circ\dac$
for all
$\a,\b,\c$.

In a similar way one can transfer to quantum categories the notions of
a covariant functor (a representation of a category), a dual category
etc. Note that in the following it will be convenient to consider both
``left'' and ``right'' coactions. Depending on it, the algebra
$\Aab$ can be interpreted as functions either on the arrows
$\a \lra \b$ or $\a\lla\b$.

\unitlength=1em
\newcommand{\ho}{\underline{\operatorname{hom}}}
\newcommand{\en}{\underline{\operatorname{end}}}
\newcommand{\gab}{g^{AB}}
\newcommand{\fkl}{f_{KL}}
\newcommand{\IV}{I_{V}}
\newcommand{\JV}{J_{V}}
\newcommand{\IW}{I_{W}}
\newcommand{\JW}{J_{W}}
\newcommand{\strelka}
{\begin{picture}(1,1)\put(1,-1){\vector(3,-2){3.3}}\end{picture}}

\section{The construction of the quantum category $\Linq$. Examples.}

Let $\Aa$ be a set of algebras. Consider the following problem.
Consider sets of algebras $\Aab$, coacting on a given set of algebras
$Aa$, in a sense that
algebra homomorphisms $\dab:\Ab\lra\Aa\otimes\Aab$ are given
for all indices , and let us look for the universal (initial)
object among such coactions.
In other words for algebras $\tilde\Aab$ coacting in the same sense on  $\Aa$,
there is a unique set of homomorphisms $\Aab\lra\tilde\Aab$ for which
the following diagram

\begin{equation}                                                              
              \begin{CD}
   \Ab\strelka@>>>\Aa\otimes\Aab\\
    @.               @VVV\\
    {}   @.       \Aa\otimes\tilde\Aab
\end{CD}
\end{equation}
is commutative.

\medskip
{\bf Theorem 2.1.} {\it The universal set of algebras $\Aab$ (if exists)
is a bialgebra. The comultiplication and counit are canonically defined
by the compatibility with the coaction.}

{\bf Proof.} The iteration of the coaction is always defined,
hence for any three indices there exists a homomorphism
$\Ac \lra\Ab\otimes\Abc\lra\Aa\otimes\Aab\otimes\Abc$.
>From the universality we obtain  homomorphisms
$\Aac\lra\Aab\otimes\Abc$. Similarly, the isomorphism
$\Aa\lra\Aa\otimes R$ yields a counit homomorphisms
$\Aaa\lra R$. The coassociativity and the counit property are easily obtained
from the uniqueness part of the universality condition.

\medskip
It is quite obvious that we can similarly treat
two sets of algebras
$\Aa,\Ba$ (or any number of them) demanding that $\Aab$
should coact on both. Or one can put additional constraints on the exact
outlook of the coactions. Then the similar universal coacting bialgebra
(if exists) will be
subject to more tight restrictions. The claim of the theorem still holds.

We shall apply this construction as follows.
Consider a set of
$\Ztwo$-graded vector spaces (over a field $k$).
In parallel to each space $V$ we consider the space $V\P$, where
$\P$ stands for the parity reversion functor. Consider the dual spaces
(the spaces of linear functions). A canonical even isomorphism
\begin{equation}
     V'\otimes V'\iso \P V' \otimes \P V', \label{p}
\end{equation}
takes the basis tensors $e^A\otimes e^B$ to $(e^A\otimes
e^B)^\P :=(-1)^{\At}(\P e^A)\otimes(\P e^B)$. Let us fix a decompositipon of
$V'\otimes V'$ to two complementary subspaces
: $V'\otimes V'=I\oplus J$.
Consider the quadratic algebras
$(V',I):=T(V')/(I)$ ... $(\P V',J^\P):=T(\P V')/(J^\P)$.
(Quotient by the ideals generated by $I$ and $J$. By $^\P$ we denote
the isomorphism  (2).) We call $(V',I)$ and $(\P V',J^\P)$ the
(polynomial) {\it function algebras on quantum superspaces }
$\Vq$ and $\VqP$.
Thus the definition of quantum superspace depends on our choice of the
decomposition
$I\oplus J=V'\otimes V'$, and the spaces  $\Vq$ and  $\VqP$
are defined simultaneously, not independently.
Denote by $x^A$ and $\x^A$
the image of basic linear functions $e^A$ and $\P e^A$ in the algebras
$\VI$ ... $\VPJ$ respectively.
We shall call them the {\it coordinates} on quantum superspaces
$\Vq$ ... $\VqP$.

\medskip
{\bf Example 2.1.} Take $I$ spanned by
$e^A\otimes e^B-(-1)^{\At\Bt}e^B\otimes e^A$
(the basis of the skew-symmetric tensors), and take $J$ spanned by
$e^A\otimes e^B+(-1)^{\At\Bt}e^B\otimes e^A$
(the basis of the symmetric tensors). Then the relations in  $\VI$ and $\VPJ$
are plain commutativity relations: $x^A x^B=(-1)^{\At \Bt}x^B x^A$,
$\x^A\x^B=(-1)^{(\At+1)(\Bt+1)}\x^B\x^A$ (can be checked).
Thus  quantum superspaces actually include  classic ones.

\medskip
Do this for each $V$. We obtain a set of quantum superspaces
$\Vq$, $\VqP$, i.e. a set of algebras $\VI, \VPJ$
which are parametrized by triples: $\VIJ$. Denote $\Aa=\VI, \Ba=\VPJ$, where
$\a=\VIJ$, and consider a coaction
\begin{equation}
\begin{aligned}
      \d:\Ab\lra\Aa\otimes\Mab,\\
      \d:\Bb\lra\Ba\otimes\Mab.
\end{aligned}
\end{equation}
The algebras $\VI$ and $\VPJ$ inherit the $\mathbb Z$-grading from the tensor
algebra.
(Not to be confused with the parity.)  We demand that the coaction must
preserve this grading (``the linearity condition''). Then
$\d(y^K)=x^A\otimes\tak$,
$\d(\h^K)=\x^A\otimes\sak $, where $x^A, y^K, \x^A, \h^K$ are coordinates on
$\Vq, W_q, \VqP, W_q\P$ respectively, and $\tak, \sak$ are certain elements
of the algebra
 $\Mab$. We also demand that $\tak=\sak$ (``compatibily with the functor
$\P$'').

{
\renewcommand{\dab}{\d_A{}^B}
\renewcommand{\dac}{\d_A{}^C}
\newcommand{\dbd}{\d_B{}^D}
\medskip
{\bf Theorem 2.2.} {\it
For the set of pairs
$\VI, \VPJ$, with the given constraints (linearity and commuting with $\P$),
there  exists the universal coacting.
This is a bialgebra $M=(\Mab)$, where all  $\Mab$ are quadratic algebras.
The choice of bases in $V$ and $W$ determines the generators
$\tak\in\Mab$, where $\a=(V,\IV,\JV), \/\b=(W,\IW,\JW)$, and
\begin{equation}
          \d(y^K)=x^A\otimes\tak, \quad   \d(\h^K)=\x^A\otimes\tak.
\end{equation}
The quadratic relations in $\Mab$ are the following:

\begin{equation}
  \begin{aligned}
        (-1)^{\Bt\Kt}f^{AB}_{(J)}\fkl^{(I)}\tak\tbl=0,\\
        (-1)^{\Bt\Kt}f^{AB}_{(I)}\fkl^{(J)}\tak\tbl=0, \label{2.5}
  \end{aligned}
\end{equation}
where we denote by $f^{(I)}, f^{(J)}$ the basis tensors in $I\oplus
J=V'\otimes
V'$ (and the same for $W$) and by $f_{(I)}, f_{(J)}$ the corresponding
elements of the dual basis.
The comultiplication
$\D:\Mac\lra\Mab\otimes\Mbc$
and the counit $\e:\Maa\lra k$ are defined by the standard formulas
\begin{align}
        \D(\tas)&=\tak\otimes\tks,\\
        \e(\tab)&=\dab
\end{align}
}

{\bf Proof.}
{
\newcommand{\bard}{\bar\delta}
We shall show that the relation~(\ref{2.5}) is valid in
any coacting
$\Aab$. Any coaction is defined by the formula~(\ref{4}), with some $\tak$.
That $\delta$ is homomorphic is equivalent to the condition
$\bard(\IW)\subset\IV\otimes\Aab\subset(V'\otimes V')\otimes\Aab$ and
$\bard(\JW^\P)\subset\JV^{\P}\otimes\Aab\subset(\P V'\otimes\P
V')\otimes\Aab$,
where $\bard$ stands for the ``covering'' coaction on the free tensor
algebras,
which exists independently of the structure
of the algebra $\Aab$. Consider
first the conditions following from the
relation $(I)$. To avoid confusion with the elements of the tensor product
of two algebras we shall omit the symbol $\otimes$ for the tensors
on a given space
(so below $e_Ae_B:=e_A\otimes
e_B$ etc). Suppose $f=\fkl e^Ke^L\in \IW\subset W'\otimes W', \quad g=\gab
e_Ae_B\in \Ann\IV\subset V\otimes V$. Then for any $f$ and $g$ it is
necessary that
$
0=\langle g,\bard(f)\rangle=\gab\fkl\langle
e_Ae_B,(e^C\otimes\tck)(e^D\otimes\tdl)\rangle=\gab\fkl(-1)^{\Dt(\Ct+\Kt)}
\langle
e_Ae_B,e^Ce^D\otimes\tck\tdl\rangle=\gab\fkl(-1)^{\Dt(\Ct+\Kt)}(-1)^
{\Bt\Ct}\dac\dbd\tck\tdl=(-1)^{\Bt\Kt}\gab\fkl\tak\tbl
$.
Thus the condition $\bard(\IW)\subset\IV\otimes\Aab$ is equivalent to
\begin{equation}
   (-1)^{\Bt\Kt}\gab\fkl\tak\tbl=0,\label{2.8}
\end{equation}
where as $f$ the basis elements of $\IW$ can be taken and as $g$
the basis elements of $\Ann\IV$.
Now we notice that the condition $\bard(\JW^\P)\subset\JV^\P\otimes\Aab$
is equivalent to $\bard(\JW)\subset\JV\otimes\Aab$.
This follows from the fact that $\bard$
commutes with the isomorphism~(\ref{p}), i.e.
$(\bard f)^\P=\bard(f^\P)$ for any $f\in
W'\otimes W'$. This immediately follows from
the equality $\bard(\P
w')=\P\bard(w')$ (the $\P$-symmetry of the coaction).
Thus it is proven that in any coacting
the relations~(\ref{2.8}) are valid, for $f\in\IW, \quad g\in\Ann\IV$
or $f\in\JW, \quad g\in\Ann\JW$. Taking into account that $I$ and $J$ are
complementary one rewrite the relations as~(\ref{2.5}).
Now take the matrix entries $\tak$ as independent variables
and consider the associative algebras generated
by them with the defining relations~(\ref{2.8}).
The set of algebras obtained will be universal.
Indeed, for an arbitrary coacting
$\Aab$ change the notation of the matrix elements to
$\sak$, then for any homomorphism commuting with the coaction
by necessity $\tak\mapsto\sak$. But this formula
defines the homomorphism uniquely
(since $\tak$ are  generators)
and it is well-defined because of the relation proven above.
\qed
}
}

\bigskip
{\bf Remarks.} {\bf 1.} As follows from the proof, the universal coacting
for the pairs of algebras
 $T(V')/(I), T(\P V')/(J^{\P})$ and for the pairs of algebras
$T(V')/(I), T(V')/(J)$ will be the same. Our preference is explained by  the
fact that in the classical case (see above) both algebras
$T(V')/(I)$ and $T(\P V')/(J^{\P})$ are commutative (in $\Ztwo$-graded sense),
while the algebra $T(V')/(J)$ is not commutative (odd generators
anticommute with even ones). The same alternative exists for the definition of
the exterior algebra, see~\cite{ber}.

\smallskip
{\bf 2.} If for each $V$ we fix a family of  complementary subspaces
$I_k^V\subset V'\otimes V'$ and look for a universal coaction for the algebras
 $T(V')/(I_k^V)$ (with the same form of the coaction on the generators,
 $\bar\d(e^K)=e^A\otimes\tak$, for all
 $k$), then we shall obtain a bialgebra $M=(\Mab)$, with the following
commutational relations for each of
$\Mab$ (matrix entries $\tak$ serve as generators):
\begin{equation}
          (-1)^{\Bt\Kt}\gab_{(k)}\fkl^{(k)}\tak\tbl=0,
\end{equation}
$f^{(k)}$ belongs to the basis of $I_k^W$, $g_{(k)}$ to the basis of
$\Ann I_k^V$. Here the
``numbers'' of objects are $\a=(V,I_1^V,\dots ,I_s^V)$ and
$\b=(W,I_1^W,\dots ,I_s^W)$.

\smallskip
{\bf 3.} The relations (\ref{2.5}) are linear independent.
It easy to calculate that for the classical values of
$\dim I_V$, $\dim J_V$, $\dim I_W$, $\dim J_V$
the dimension of the quadratic part of the ``quantum algebra'' $\Mab$
will be classical too. The question whether the same will be true for
the higher order terms is highly non-trivial. We discuss it in the next
section.

\smallskip
{\bf 4.} In  quantum group theory the universal coaction method
was proposed by Yu.A.Kobyzev
(see~\cite{man1}). Originally a single quadratic algebra $A=T(V')/(I)$
was considered. This gives only ``half'' of the necessary commutational
relations on quantum (semi)group. This difficulty has been resolved in a
rather artificial manner by using the ``dual'' quadratic algebra $A^!$
(see the remark after Example 2.5 below). In particular that
implied that a single subspace $I\subset V'\otimes V'$ played the role of
the ``quantization parameter''. The application of  families of
complementary subspaces is described in~\cite{dem2}, in purely even case.
(The functor $\P$ is not used in ~\cite{dem2}.)
It is worth noting that in
~\cite{man1} a quadratic algebra
$\ho(A,B)=A^!\bullet B$ was introduced for each pair of
quadratic algebras $A$ and $B$, together with a sort of
``coproduct'' $\D:\ho(A,C) \lra \ho(B,C)\circ\ho(A,B)$ (we refer
also to ~\cite{gel} for the notation). The algebras $\ho$ are not so
interesting (they, and in particular $\en(A)=\ho(A,A)$,
indeed lack half of the relations, so they are strongly non-commutative).
But having in mind the natural homomorphism $A\circ B\lra A\otimes B$
~\cite{gel},~\cite{man1} they can be considered as an example of a quantum
category in our sense, or. more precisely, an examlple of the generalized
bialgebra (in the sense of Definition 1.2).
The interpetation given in ~\cite{man1} was completely different.
It was based on the ``internal $\Hom$'' formalism of the ``rigid tensor
categories'', which are just actual categories endowed with the additional
structure.

\bigskip
{\bf Definition 2.1.} The bialgebra $M=(\Mab)$, defined in  Theorem 2.2,
will be called the {\it algebra of {\em(polynomial)} functions on the
quantum category}
$\Linq$.

\bigskip
{\bf Example 2.2.} Consider classical (not deformed)
spaces. For every $V$ the subspaces $\Ann I, \Ann J \subset
V\otimes V $ are spanned by the tensors $e^A\otimes
e^B+(-1)^{\At\Bt}e^B\otimes
e^A$ and $e^A\otimes e^B-(-1)^{\At\Bt}e^B\otimes e^A$. By Theorem 2.2
we obtain the following relations for the matrix elements of the
homomorphism from $V$ to $W$:
{
\newcommand{\daa}{\d_A{}^{A'}}
\renewcommand{\dab}{\d_A{}^{B'}}
\newcommand{\dbb}{\d_B{}^{B'}}
\newcommand{\dba}{\d_B{}^{A'}}
\newcommand{\dkk}{\d_{K'}{}^{K}}
\newcommand{\dll}{\d_{L'}{}^L}
\newcommand{\dkl}{\d_{K'}{}^L}
\newcommand{\dlk}{\d_{L'}{}^K}
\renewcommand{\tak}{t_{A'}{}^{K'}}
\renewcommand{\tbl}{t_{B'}{}^{L'}}
\newcommand{\bt}{\tilde {B'}}
\newcommand{\kt}{\tilde {K'}}

\begin{eqnarray*}
   (\daa\dbb+(-1)^{\At\Bt}\dba\dab)(\dkk\dll-(-1)^{\Kt\Lt}\dkl\dlk)
                                             (-1)^{\bt\kt}\tak\tbl=0, \\
   (\daa\dbb-(-1)^{\At\Bt}\dba\dab)(\dkk\dll+(-1)^{\Kt\Lt}\dkl\dlk)
                                             (-1)^{\bt\kt}\tak\tbl=0.
\end{eqnarray*}
}
Then, after summation, we obtain:
 $$
(-1)^{\Bt\Kt}\tak\tbl-(-1)^{\Kt\Lt+\Bt\Lt}\tal\tbk+(-1)^{\At\Bt+\At\Kt}\tbk
\tal-(-1)^{\At\Bt+\Kt\Lt+\At\Lt}\tbl\tak=0
 $$
and a similar equality  with the opposite signs
before the second and the third terms.
Summing and substracting these equalities, we get
 $$
        \tak\tbl-(-1)^{(\At+\Kt)(\Bt+\Lt)}\tbl\tak=0,
 $$
for any $A,B,K,L$, i.e. the usual commutativity condition.

\medskip
{\bf Example 2.3.} Take one-dimensional  $k$ as  $W$ with
 $0$ and $k\iso k\otimes k$
as $I_W, J_W$ respectively. Then by Theorem 2.2 we get
\begin{equation}
   \gab_{(J)}t_A t_B=0
\end{equation}
 as commutational relations for the coefficients of even {\it linear
forms} on $\Vq, \VqP$. Here $g_{(J)}\in \Ann J_V$. In the same way, for the
coefficients of {\it odd linear forms} $\theta_A, \quad \tilde
\theta_A=\At+1$, we obtain the relations
\begin{equation}
   (-1)^{\Bt}\gab_{(I)}\theta_A\theta_B=0,
\end{equation} where $g_{(I)}\in \Ann I_V$.
In other words, if one identifies
a quantum space with a triple $(V,I,J)$, then its  {\it quantum dual space
} will be $(V',\Ann J, \Ann I)$. Pairings: $\langle x,t
\rangle=x^A\otimes t_A, \quad \langle x, \theta \rangle =x^A\otimes\theta_A$.

\medskip
{\bf Example 2.4.} From the previous example one can deduce the
commutational relations for the coefficients of
(even) {\it bilinear forms},
considered as homomorphisms to dual space (``lowering indices''):
\begin{equation}
\begin{aligned}
       (-1)^{\Bt\Ct}g_{(I)}^{AC}g_{(J)}^{BD}t_{AB}t_{CD}=0,\\
       (-1)^{\Bt\Ct}g_{(J)}^{AC}g_{(I)}^{BD}t_{AB}t_{CD}=0,
\end{aligned}
\end{equation}
where $\tilde t_{AB}=\At +\Bt, \quad g_{(I)}, g_{(J)}$ are tensors in $\Ann I,
\Ann J \subset V\otimes V$ respectively. The bilinear forms are:
$\langle x_1\mid T\mid x_2\rangle = x_1^A\otimes t_{AB}\otimes x_2^B$
and the similar for $\x_i$.

\medskip
{\bf Example 2.5.} Consider the quantum spaces with the following
commutational relations:
\begin{equation}
\begin{aligned}
    x^A x^B-q^{AB}x^B x^A&=0,\\
    \x^A\x^B+(-1)^{\At+\Bt}p^{AB}\x^B\x^A&=0.
\end{aligned}      \label{sudb}
\end{equation}
Here $\qab=(\qba)^{-1}$, $\pab=(\pba)^{-1}$, $\qaa=\paa=(-1)^{\At}$.
(In the ``classical limit'' $\qab,\pab\lra(-1)^{\At\Bt}$.)
Look for the relations in the function algebra
for the corresponding ``full subcategory''
of
$\Linq$. Here the objects are  triples
$(V,\pv,\qv)$,
where
{\tolerance=500
${\pv=(\pab)=(\pab_{(V)})}$, ${\qv=(\qab)=(\qab_{(V)})}$ are the matrices
of the parameters. The subspace $I$ is spanned by the tensors
${e^A\otimes e^B-\qab e^B\otimes e^A}$,  $J$ by the tensors
${e^A\otimes e^B+\pab e^B\otimes e^A}$,  $\Ann I$  by the tensors
${e_A\otimes e_B+\qban e_B\otimes e_A}$, and $\Ann J$ by the tensors
${e_A\otimes e_B-\pban e_B\otimes e_A}$.
As not to contradict the tensor notation, we have introduced here
the parameters with the lower indices:
$\qabn:=\qab, \quad \pabn:=\pab$.
The complementarity of $I$ and $J$
is equivalent to
\begin{equation}
      \qab+\pab\neq0, \label{comp}
\end{equation}
(for any $V$). Similarly to Example 2.1 we obtain two relations
for the matrix elements of a
``linear map of the quantum space
 $(V,\pv,\qv)$ to the quantum space
$(W,\pw,\qw)$'':
 $$
   \begin{aligned}
(-1)^{\Bt\Kt}\tak\tbl-\qkl(-1)^{\Bt\Lt}\tal\tbk+\qban((-1)^{\At\Kt}\tbk
\tal-\qkl(-1)^{\At\Lt}\tbl\tak)=0\\
(-1)^{\Bt\Kt}\tak\tbl+\pkl(-1)^{\Bt\Lt}\tal\tbk-\pban((-1)^{\At\Kt}\tbk
\tal+\pkl(-1)^{\At\Lt}\tbl\tak)=0
    \notag
   \end{aligned}
 $$
Here $\qv=(\qab), \pv=(\pab), \qw=(\qkl), \pw=(\pkl)$.
Provided \ref{comp}, these defining relations can be identically transformed
to the following final form:
\begin{multline}
  \tak\tbl-\frac{\pban+\qban}{\plk+\qlk}(-1)^{\At\Lt+\Bt\Kt}\tbl\tak=\\
\frac{\pban\plk-\qban\qlk}{\plk+\qlk}(-1)^{(\At+\Bt)\Kt}\tbk\tal  \label{rel}
\end{multline}
for any $A, B, K, L$. Notice that for the elements of one column or of
one row
these relations reduce to
\begin{align}
    \tak\tbk-\frac{\pban(1+(-1)^{\Kt})+\qban(1-(-1)^{\Kt})}{2(-1)^
                           {(\At+\Bt+1)\Kt}}\tbk\tak=0,  \label{rel1}\\
    \tak\tal-\frac{2(-1)^{(\At)(\Kt+\Lt+1)}}{\plk(1-(-1)^{\At})+\qlk
                         (1+(-1)^{\At})}\tal\tak=0.       \label{rel2}
\end{align}

\bigskip
{\bf Remarks. 1.} The particular case of (\ref{rel}) are the relations in
the function algebra on the
``multiparameter deformation'' of the general linear supergroup.
The original relations of the form  (\ref{sudb}) for quantum linear spaces
(in the even case) together with the relevant quantization of $\GLn$
is due to Sudbery~\cite{sud} (see the surveys~\cite{dem2}, ~\cite{dem3}).
Before Manin~\cite{man2} introduced a multiparameter quantization
for supergroup
$\GLnm$ starting from the relations similar to  (\ref{sudb})
with  $(\qab)^{-1}$ instead of $\pab$ (in our notation).
There was a certain confusion in the papers~\cite{man1},
~\cite{man2}, ~\cite{dem1} between a vector space and its dual
and as a consequence the analogue of the second equation of (\ref{sudb})
was associated with the ``dual quadratic algebra''
$A^!$ while the coaction was taken as if it had been for the initial space.

The structure of the formulas (\ref{rel}-\ref{rel2})  is similar
to the relations for the standard (single-parameter) deformation
of the group $\operatorname{GL}\,(2)$ ~\cite{rtf}. The commutational relations
for the multiparameter deformations of $\GLnm$ are commonly presented
in a more cumbersome form
~\cite{man1}, ~\cite{man2},
~\cite{dem1}, ~\cite{dem2},~\cite{dem3}.

{\bf 2.} In the paper~\cite{dem1} Demidov defined for two spaces $V$ and $W$
an algebra close to that defined above,
for $\pab=(\qab)^{-1},
\quad \pkl=(\qkl)^{-1}$, see the previous remark. But he never considered
a comultiplication
except for the case of a single space with fixed parameters
$\qab$ (i.e. for ``quantum semigroups'').

\bigskip
{\bf Example 2.6.} For quantum space with the relations (\ref{sudb})
the commutational relations for the dual space (see Example 2.3)
are:
\begin{equation}
   \begin{aligned}
       t_A t_B-(\pabn)^{-1}t_B t_A=0,\\
       \theta_A \theta_B-(-1)^{\At+\Bt}(\qabn)^{-1}\theta_B\theta_A=0.
   \end{aligned}
\end{equation}

\medskip
{\bf Example 2.7.}  Similarly for  even bilinear forms (see Example 2.4):
{
\renewcommand{\tab}{t_{AB}}
\newcommand{\tad}{t_{AD}}
\newcommand{\tcd}{t_{CD}}
\newcommand{\tcb}{t_{CB}}
\begin{multline}
  \tab\tcd-\frac{\pcan+\qcan}{\qbdn+\pbdn}(-1)^{\At\Dt+\Bt\Ct}\tcd\tab=\\
\frac{\pcan\qbdn-\qcan\pbdn}{\qbdn+\pbdn}(-1)^{(\At+\Ct)\Bt}\tcb\tad
\end{multline}
}

\bigskip
For convenience we shall also write down
the relations (\ref{sudb}--\ref{rel2}) for the subcategory of purely
even spaces
($V$ is purely even, $V\P$ is purely odd):
{
\renewcommand{\qlk}{q^{lk}}
\renewcommand{\plk}{p^{lk}}
\renewcommand{\qba}{q_{ba}}
\renewcommand{\pba}{p_{ba}}
\renewcommand{\tbl}{t_b{}^l}
\renewcommand{\tak}{t_a{}^k}
\renewcommand{\tbk}{t_b{}^k}
\renewcommand{\tal}{t_a{}^l}
\begin{equation}
    \begin{aligned}
        x^a x^b-\qab x^b x^a&=0,\\
        \x^a\x^b+\pab \x^b \x^a&=0,  \label{evsudb}
    \end{aligned}
\end{equation}
\begin{align}
   \tak\tbl-\frac{\pba+\qba}{\plk+\qlk}\tbl\tak&=\frac{\pba\plk-\qba\qlk}
              {\plk+\qlk}\tbk\tal,\\
 \tak\tbk-\pba\tbk\tak=&0,\\
 \tak\tal-(\qlk)^{-1}\tal\tak&=0,
\end{align}
}

\bigskip
Definition 2.1 makes sense for just  a single vector space
$V$. Here we obtain
the ``full subcategory'' of $\Linq$ consisting of all
linear transformations between
various quantum deformations of $V$ defined by the decomposition $I\oplus
J=V'\otimes V'$.
This quantum category is the  adequate quantization of the semigroup $\End V$.
We see that the quantum categories language is more to the point
here. It is more flexible than the language of quantum (semi)groups.

\medskip
{\bf Example 2.8.} Fix a space $V,\quad \dim V=2$ and let us change
the quantization parameters in (\ref{evsudb}). Suppose $p^{21}=p$,
$q^{21}=q$,$p=\pa$, $q=\qa$. We arrive to a quantum category,
with indices $\a, \b$ that number the parameters $p$ and $q$ as objects.
The commutational relations for the matrix elements
of a ``homomorphism from $(V,\a)$ to $(V,\b)$'' in the  natural notation
are:
\begin{align}
   ab-\qb^{-1}ba&=0, \label{ab} \\
   ac-\pa ca&=0,       \label{ac} \\
   ad-\frac{\pa+\qa}{\pb+\qb}da&=\frac{\pa\pb-\qa\qb}{\pb+\qb}cb,
     \label{ad} \\
   bc-\frac{\pa+\qa}{\pb+\qb}\pb\qb cb&=\frac{\pa\qb-\qa\pb}{\pb+\qb}da,
        \label{bc}\\
   bd-\pa db&=0, \label{bd}\\
   cd-\qb^{-1} dc&=0.  \label{cd}
\end{align}
From~(\ref{ad}) and~(\ref{bc}) the relation between $b$ and $c$ can be
expressed also as
\begin{equation}
      bc-\frac{\pb+\qb}{\pa+\qa}\pa\qa
cb=\frac{\pa\qb-\qa\pb}{\pa+\qa}ad.
\end{equation}
For $\a=\b$ these relations reduce to the commutational relations on
the quantum group $\operatorname{GL}(2)$.

\medskip
{\bf Example 2.9} ({\sl the determinant on  a quantum category}).
In the situation of the previous example
consider the ``area form'' $\x^1\x^2$.
Substituting
$\x_{\b}^k=\x_{\a}^a t_a{}^k$, we obtain
\begin{equation}
  \d(\x_{\b}^1\x_{\b}^2)=\x_{\a}^1\x_{\a}^2 \otimes \detab(T).
\end{equation}
The factor denoted by $\detab(T)$ is called, by definition,
the {\it determinant} of the matrix $T=(t_a{}^k)$.
Here $\x_{\a}^a $ are the coordinates on $(\P V,\a)$ and $ \x_{\b}^k$
the coordinates
on $(\P V,\b) $. From the definition and the
commutational relations
we obtain the equivalent formulas:
\begin{equation}
\begin{split}
  \detab(T)&= ad-\pa cb\\
           &=\frac{\pa+\qa}{\pb+\qb}(da-\qb cb)      \\
           &=\frac{1+\pa\qa^{-1}}{\pb+\qb}(\qb ad-bc) \\
           &=\pb^{-1}(\pa da-bc).
\end{split}
\end{equation}
For any ``morphisms $\a\lra\b, \quad \b\lra\c$''
the identity
\begin{equation}
   \detac(\Tac)=\detab(\Tab)\detbc(\Tbc),
\end{equation}
holds.  Here
$\Tac=\Tab\Tbc$. Thus
$\det$ defines a one-dimensional representation of the quantum category in
consideration.
Mind that the function $\detab$ is defined not uniquely but up to
the ``adding of a coboundary'', i.e. up to a factor
$f_\a f_\b^{-1}$, where $f_\a$ is a non-vanishing scalar function
of the parameters
$\pa, \qa$. (Due to the non-uniqueness of the basis form
 $\x_{\a}^1\x_{\a}^2$). The subcategory with $\det$
cohomologous to zero is a correct
quantum analogue of the
group
$\operatorname{SL}$.

\newcommand{\cab}{c_{AB}}
\newcommand{\ckl}{c_{KL}}
\newcommand{\eab}{\e_{AB}}
\newcommand{\eba}{\e_{BA}}
\newcommand{\ebc}{\e_{BC}}
\newcommand{\eac}{\e_{AC}}
\newcommand{\ekl}{\e_{KL}}
\newcommand{\Bplus}{B_+}
\newcommand{\Bminus}{B_-}
\newcommand{\Bnmkl}{B_{NM}{}^{KL}}
\newcommand{\Babcd}{B_{AB}{}^{CD}}
\newcommand{\Bod}{B^{12}}  
\newcommand{\Bdt}{B^{23}}  
\newcommand{\la}{\l_{\a}}
\newcommand{\lb}{\l_{\b}}
\newcommand{\lk}{\l_{k}}
\newcommand{\lj}{\l_{j}}
\newcommand{\sba}{\operatorname{sign}(B-A)}
\newcommand{\sab}{\operatorname{sign}(A-B)}
\newcommand{\slk}{\operatorname{sign}(L-K)}
\newcommand{\Ik}{I_k}
\newcommand{\Il}{I_l}
\newcommand{\IkV}{I_k^V}
\newcommand{\IkW}{I_k^W}
\newcommand{\PkW}{P_k^W}
\newcommand{\PkV}{P_k^V}
\newcommand{\PjW}{P_j^W}
\newcommand{\PiV}{P_i^V}
\newcommand{\Pk}{P_k}

\section{ The Poincare -- Birkhoff -- Witt property.
\mbox{$R$-matrix} formulation and the connection with the Yang -- Baxter
equation.}

Consider the quantum category $\Linq$ with triples
$(V,I.J)$, $I\oplus J={V'\otimes V'}$ as objects.
It is possible to consider even more general case
of the quantum category with objects like
${(V,I_1,\dots,I_s)}, \quad \bigoplus \Ik={V'\otimes V'}$,
the number $s$ is fixed (see the remark after Theorem 2.2).
Denote it by $\Linqs, \, s\geqslant 2$. Then
$\Linq=\Linq^{(2)}$.

\medskip
{\bf Theorem 3.1.} \  {\it Let the matrix entries of  $\Tab=(\tak)$
generate the algebra
$\Mab$, which is the function algebra on homomorphisms from
$\a={(V,I_1^V,\dots,I_s^V)}$ to
$\b={(W,I_1^W,\dots,J_s^W)}$, in the quantum category $\Linqs$. Let
$\Tab^1:=\Tab\otimes1,\quad \Tab^2:=1\otimes\Tab$.
Then the defining commutational relations for functions on
$\Linqs$ can be presented in the following
``$R$-matrix form'':
\begin{equation}
    \Ba(\Tab^1\Tab^2)=(\Tab^1\Tab^2)\Bb, \label{1}
\end{equation}
where for each $\a$ the matrix $\Ba$
is the linear combination of the projectors
on the subspaces $\Ik$ along the sum
$\bigoplus\limits_{l\neq k}^{}\Il$, where the coefficients  $\lk$
are pairwise
different and independent on
$\a$. }

\medskip
{
\renewcommand{\d}{\bar\delta}
\newcommand{\djj}{\d_{jj}}
\newcommand{\dkj}{\d_{kj}}
{\bf Proof.} Consider the partitions of unity $1=\sum \PkV$, $1=\sum \PkW$
by the projectors
 corresponding to the direct sum decompositions
$V'\otimes V'=\bigoplus\IkV$, $ W'\otimes W'=\bigoplus\IkW$. Let
$A_{\a,k}=T(V)/(\IkV)$,  $A_{\b,k}=T(W)/(\IkW) $. The homorphisms  $\delta:
A_{\b,k}\lra A_{\a,k}\otimes\Mab, \quad k=1,\dots,s$, (the coaction)
are covered by the map of the tensor algebras
$\d$, and the commutational relations in $\Mab$ are defined by the condition
$\d(\IkW)\subset \IkV\otimes\Mab $ for all $k$. Any linear map
$A:W'\otimes W'\lra V'\otimes V'\otimes\Mab$ can be uniquely decomposed as
$A=\sum\PiV\! A\PjW=\sum A_{ij}$. The condition above is equivalent
to that the operator
$\d$ is ``diagonal'' on the tensor square:
$\d=\sum\djj$. Let us introduce
$\Ba=\sum\lk\PkV, \Bb=\sum\lk\PkW$, where $\lk\neq\lj, \, k\neq j$. Then
the difference of the left- and right-hand sides of
~(\ref{1})
is
$\Bb\circ\d-\d\circ\Ba=\sum\limits_{k,j}(\lk-\lj)\dkj$, which vanishes if
and only if
$\dkj=0$ for $k\neq j$. That is just
 the condition defining the commutational relations in
 $\Mab$. \qed
}

\bigskip
In the index notation the relation (\ref{1})
can be rewritten as follows:
\begin{equation}
       \tan\tbm\Bnmkl=\Babcd\tck\tdl,    \label{2}
\end{equation}
where  $\Ba=(\Babcd),\quad  \Bb=(\Bnmkl)$.

The claim of the theorem is analoguous to the
corresponding result for quantum groups,
where the relations are determined by
 a single matrix $B$.

\medskip
For any quantum space $(V,I_1,\dots,I_s)$ the projectors $\Pk$
can be expressed from the matrix $\Ba$ (``$R$-matrix'')
as the projections on its eigenspaces.
The commutational relations in the algebras
$A_k=T(V')/(I_k)$ are also expressed in terms of the operator $\Ba$.
Suppose one has the commutational relations in the $R$-matrix form.
It is common then to relate ``the Poincare -- Birkhoff -- Witt property''
(the equality of the graded dimension of function algebra on quantum space
or quantum group to that of the corresponding classical polynomial algebra)
with the identity
\begin{equation}
  \Bod\Bdt\Bod=\Bdt\Bod\Bdt \label{YB}
\end{equation}
for the matrix $B$, where $\Bod=B\otimes1, \Bdt=1\otimes B$
(the Yang -- Baxter equation). Strictly speaking, no general statement exists
here and one can only say that a certain ``relation'' takes place.
(As Joseph Donin told me after this work was completed,
a sort of theory can be actually developped dealing with the quantum
algebra dimension in connection with the YB equation. This is the topic
of the recent paper~\cite{don} by him and Steve Shnider.)

Consider the case of $\Linq=\Linq^{(2)}$, $I_1=I, I_2=J$.
We shall normalize the matrices
$B$ so as $B=P_1-\l P_2,\quad \l\neq-1$, and we shall look for $\l$ suiting
(\ref{YB}).
In general this is an overdetermined problem.
One can notice
(see, for example, ~\cite{dem2}), that if the solutions exist there
may be either exactly two of them:
$\Bplus=P_1-\l P_2$ and $\Bminus=P_1-\l^{-1}
P_2$, where $\l\neq1$, or just one: $B=P_1-P_2$, in the case the equation
(\ref{YB}) is satisfied for
$\l=1$. (In the last case an additional identity holds, $B^2=1$.)
Applying this to our situation, we obtain that if for two quantum spaces
$\a=(V,I_1^V,I_2^V)$ and $\b=(W,I_1^W,I_2^W)$ the normalized
matrices $\Ba$ and $\Bb$ suit the Yang -- Baxter equation, and we fix the
choice of $\Ba$ and $\Bb$,
then the commutational relations in $\Mab$
can be presented in the ``$R$-matrix form''(\ref{1}),(\ref{2}) with the
given $\Ba$ and $ \Bb$ if and only if
$\la=\lb$ or $\la=\lb^{-1}$.

\medskip
Let us turn now to problem of the graded dimension of the algebras $\Mab$.
We consider the subcategory with defining commutational relations for
quantum spaces of the form (\ref{sudb}).
Let $\a=(V,P_V,Q_V),
\quad \b=(W,P_W,Q_W)$ with the matrices of parameters $P_V=(\pab), Q_V=(\qab),
P_W=(\pkl), Q_W=(\qkl)$. For arbitrary orderings of the tensor indices in
$V$ and $W$ let $\sba$ equal $+1, 0, -1$ if $B>A,\
B=A,\ B<A$ respectively. In a similar way we define the function $\slk$.
The matrix entries will be ordered
by rows: $\tak<\tbl$, if $A<B$ or if $A=B,\quad K<L$. The {\it
ordered monomial\/} is defined as usual to be such that
the sequence of letters is not decreasing and any odd variable can appear
no more than once.

\medskip
{\bf Theorem 3.2.} \ {\it The algebra $\Mab$ is of classical dimension iff
the bases in $V$ and $W$ can be ordered in such a way that
the following equations hold
\begin{align}
   \pab&=\qab c_\a^{\sba},\label{a} \\
   \pkl&=\qkl c_\b^{\slk},
   \label{4}
\end{align}
for some constants $c_\a,\/c_\b $, not equal to $0$ or $-1$,  and either
\begin{equation}
     c_\a=c_\b,      \label{eq}
\end{equation}
or
\begin{equation}
     c_\a=c_\b^{-1}, \label{inv}
\end{equation}
In this case the ordered monomials in matrix entries span the algebra $\Mab$.}

\medskip
{\bf Proof.} It actually consists of two parts.

First choose an arbitrary ordering of the tensor indices.
We shall show that the ordered monomials span the whole algebra
$\Mab$
for any values of the parameters $\pab, \qab, \pkl,\qkl$.
Indeed, consider an arbitrary monomial in
$\tak$.
Consider the elements of the first row in it.
Let $\t$ stand for the number of the inversions with the elements
of other rows. That is the number of cases when the elements of the first
row appear to the right of the elements of other rows.
Take any neighbouring pair of elements
of the sort
${\tak t_1{}^L}, A>1 $.
By the commutational relations (\ref{rel})
it is possible to
substitute it by a linear combination
of the products  $t_1{}^L \tak$ ... $t_1{}^K\tal$,
and the number $\t$ will decrease by one.
Thus a finite iteration of this procedure
leads to a linear combination of the monomials
such that in each of them the elements of the first row
stand to the left of the product of the elements of the other rows.
We can apply the same method to these products, now taking the second row
instead of the first, etc.
(If at some step the elements of a certain row are absent, then we take the
next row.) As there are finite number of rows, after a finite steps we shall
obtain a linear combination of monomials with the correct order of the row
numbers. The elements inside the same row (they all are neighbours now)
are put in the correct order by the formulas (\ref{rel2}). Now the arbitrary
monomial from which we have started is expressed as linear combination of
the ordered mononials.

Such algorithm for a quantum deformation
of the supergroup $\GLnm$ was described by Manin
~\cite{man2} with the reference  to Yu.Kobyzev.

The second part of the proof deals with the linear independence of
the ordered monomials. This  is much
more non-trivial problem. The so-called ``diamond lemma'' of the combinatorial
theory of rings
(see ~\cite{lat}, ~\cite{dem3}) implies that for quadratic algebras it is
necessary and sufficient to check the independence only
for the monomials of the third order.
(The lemma provides in this case the linear independence also for all
ordered monomials of the degree~ $\geqslant 3$.)
Let us look for the criterion of the independence for ordered cubic
monomials in
$\tak$.

The linear independence of ordered monomials is equivalent to the uniqueness
of the ``normal form'' of any arbitrary monomial, i.e. to the independence
``on path'' of the result of its expression
as a combination of  the ordered monomials.
We need to finger all cubic monomials
with the broken order and for each of them to compare all the possible ways
to put them in the normal order.
Let ${a<b<c}$ be the  arbitrary letters symbolizing
the matrix entries.
Then the cubic monomials with the broken order are the following:
$acb$, $bac$,
$bca$,
$cab$, $cba$ (three letters different) and $aba$, $ba^2 $ (two letters
different).
In the process of putting them into order
the new letters can arise.
Which ones (together with the ways of putting into order) actually
depends on the position of the elements
$a,b,c$ in the matrix $(\tak)$. It is most convenient
to draw pictures on which the matrix entries
are shown as the vertices of an orthogonal lattice.
By (\ref{rel}) the ``commutator closure'' of the given set of letters
consists of the vertices of the least sublattice containing this set.
Other letters cannot appear in the process of putting the monomial into
normal form.
For three pairwise different entries we obtain one picture with
six additional letters (if the initial elements stand on three different
rows) and two pictures with three additional letters (if initially we have
elements in two different rows).
For two different letters a picture with two additional letters appears.
(The pictures with the initial letters in one row do not include any
additional letters and they cannot produce any ambiguity.
The staightforward checking of all possible ways leading
to the sum of ordered monomials in each of these
$5\times 3+2=17$ cases shows that the non-uniqueness of the normal form
is apriori possible only for the monomials like
 $cba$ (for all three variants of the position of the initial elements),
i.e. when branching of the algorithm takes place at  the very first step.
(Have in mind that the conclusion is strongly based on the exact form of
the commutational relation.)
Happily no branching is possible here at the following steps.
In the most complex case with six additional letters two branches of the
algorithm look as follows:
$cba\lra{(bca,xya)}\lra{(bac,buv,xay,xtv)}\lra
{(abc,tsc,ubv,txv,axy,usy,txv,ubv)}$ and
$cba\lra{(cab,cts)}\lra{(acb,uvb,tcs,uys)}\lra
{(abc,axy,ubv,usy,tsc,txv,usy,ubv)} $.
Here all the arising monomials are successively written down.
In two other cases the branches are a bit shorter.
The final list of the ordered monomials is the same for both branches.
That means that the uniqueness of the normal form is equivalent
to the agreement of numeric factors in the similar terms.
The explicit calculation leads us to the  12 equations
for the parameters $\pab,\qab,\pkl,\qkl$. Five of these equations happen to
be identities.
Let us give as example two of the non-trivial equations:
\begin{multline}
    \frac{(\pabn\qln-\qabn\pln)(\pacn\pkn-\qacn\qkn)(\pbcn+\qbcn)}{(\pkn+\qkn)
(\pln+\qln)} +\\
    \frac{(\pabn+\qabn)(\pacn\pkl-\qacn\qkl)(\pbcn\pln-\qbcn\qln)}{(\pkl+\qkl)
(\pln+\qln)}\\=
    \frac{(\pabn\pkl-\qabn\qkl)(\pacn+\qacn)(\pbcn\pkn-\qbcn\qkn)}{(\pkl+\qkl)
(\pkn+\qkn)},    \label{7}
\end{multline}
if ${A<B<C}$, and
\begin{multline}
    \frac{2(\pabn+\qabn)(\pabn\pkn-\qabn\qkn)}
         {(\pkn+\qkn)(\pnl+\qnl)(\pln(1-(-1)^{\Bt})+\qln(1+(-1)^{\Bt}))}=\\
    \frac{2(\pabn+\qabn)(\pabn\pkn-\qabn\qkn)}
         {(\pkn+\qkn)(\pkl+\qkl)(\plk(1-(-1)^{\Bt})+\qlk(1+(-1)^{\Bt}))}+\\
    \frac{(\pabn\pkl-\qabn\qkl)(\pabn\pln-\qabn\qln)}
         {(\pkl+\qkl)(\pln+\qln)},    \label{8}
\end{multline}
if $A<B, L<N$. To get the solution we shall make a reduction.
Substitute  $K=N$ into ~(\ref{8}). Two terms cancel immediately
and we arrive to
 $(\pabn\pnl-\qabn\qnl)(\pabn\pln-\qabn\qln)=0$,
which is equivalent to
\begin{equation}
    \frac{\pabn}{\qabn}=\frac{\pln}{\qln}  \quad\text{ or }\quad
    \frac{\pabn}{\qabn}=\frac{\qln}{\pln}     \label{9}
\end{equation}
(recall that $\pln=(\pnl)^{-1}, \qln=(\qnl)^{-1}$). Let $\pabn/\qabn=\cab,
\pkln/\qkln=\ckl$. The equations ~(\ref{9}) are equivalent to
\begin{equation}
   \cab+(\cab)^{-1}=\ckl+(\ckl)^{-1}  \label{10}
\end{equation}
for any $A<B,\quad K<L$. That means that both the right-hand side and the
left-hand side of~(\ref{10}) do not depend on indices and are equal to some
constant
$\m\neq 0$. Thus for all
$A,B \quad \cab=c^{\eab}$, where $c: c+c^{-1}=\m$, $ \eab+\eba=0$,
$\eab=\pm 1$  if $A\neq B$, and the same for $\ckl=c^{\ekl}$. Let $c\neq
1$. If the matrices $(\eab), (\ekl)$ do not possess the ``transitivity"
property
($\eab=1, \ebc=1$ implies $\eac=1$), then the equations are not satisfied.
Indeed, substituting, for example, $\pabn=c\cdot\qabn$,
$\pbcn=c\cdot\qbcn$, $\pacn=c^{-1}\cdot\qacn$
into~(\ref{7})
and setting $K=L=N$, we obtain a contradiction (for $c\neq 1$).
It follows that if it is not true that $\eab=\sba, \ekl=\slk$
for some ordering of the bases
of  $V$ and $W$ (perhaps different from the originally chosen),
then the dimension of the algebra $\Mab$ is less than the classical value.
 (The ordered cubic monomials, w.r.t. the original ordering, will be linear
dependent and according to what was proven above they actually span the
subspace of all cubic functions.)
Assume now that this is true. We write the same equations on the
deformation parameters w.r.t. the new ordering.
To finish the proof it is
necessary to substitute in them the equations
~(\ref{4}), ~(\ref{eq}). Now the direct, rather laborous
check shows that all the seven equations are actually satisfied.
The case of $c=1$ is also contained here; it is indifferent to the choice
of ordering.
\qed

\bigskip
We actually proved that it is enough to consider the case
(\ref{eq}), for a suitable ordering. From the other point, it was known that
for a single space
$\a=(V,I_1,I_2)$ the possibility to choose the order of the basis in such a
way that the equation
(\ref{a}) holds for a certain $c_\a$ is equivalent both to the PBW property
for the algebra
$\Maa$ and to the Yang -- Baxter equation for
${B_\pm=P_1-(c_\a)^{\pm1}P_2}$ \ (see ~\cite{dem2},~\cite{dem3}).
The key point of our result is that the equality of ``quantum constants''
$c_\a=c_\b^{\pm1}$ for two spaces is exactly the criterion for the classical
value of the dimension of the algebra
 $\Mab$ (the ``Poincare -- Birkhoff -- Witt'' for $\Mab$).
This criterion is in excellent agreement
with our heuristic consideration above, which was based on the Yang --
Baxter equation,
while there is no direct logical connection.
For  $c_\a\neq1$ one can choose one of the two orderings compatible with
(\ref{a})  (direct or reverse),
i.e. to choose one of the two
``Yang--Baxter structures" $\Bplus$ or $\Bminus$ on $(V,I_1,I_2)$,
which could be viewed as a sort of ``orientation''. The statement of the
theorem takes it into account.

\smallskip
It is convenient to make a slight change of the notation.
Namely, instead of $c_\a$ we introduce $\l_\a: \/  c_\a=\l_\a^2$ and
rescale the parameters, $\qab:=\qab\l_\a^{\sab}$.
Then the commutational relations will become more symmetric.

Fix now a  number $\l\neq0,\pm i$. Define a quantum category $\Linqplus$,
with objects which are triples like $(V,Q,\pm1)$ with the following relations
\begin{equation}
  \begin{aligned}
       x^A x^B-\qab\l^{\pm\sab}x^B x^A=0,\\
       \x^A\x^B+(-1)^{\At+\Bt}\qab\l^{\mp\sab}\x^B \x^A=0,
  \end{aligned}     \label{space}
\end{equation}
(the basis in $V$ is defined uniquely up to the scaling of the basis vectors).
Then the matrix entries of the ``homomorphisms from $(V,Q_V,\e)$ to
$(W,Q_W,\h)$'' are subject to the following commutational relations
implied by
(\ref{rel}-\ref{rel2}) ... (\ref{space}). For the elements of one row, if
$K<L$:
\begin{equation}
  \begin{aligned}
    \tak\tal-\qkl\l^{-\h}\tal\tak&=0, \text{   $\At=0$,}\\
    \tak\tal+(-1)^{\Kt+\Lt}\qkl\l^{\h}\tal\tak&=0, \text{   $\At=1$}.
  \end{aligned}
\end{equation}
For the elements of one column, if $A<B$:
\begin{equation}
  \begin{aligned}
    \tak\tbk-(\qabn)^{-1}\l^{-\e}\tbk\tak&=0, \text{    $\Kt=0$,}\\
    \tak\tbk+(-1)^{\At+\Bt}(\qabn)^{-1}\l^{\e}\tbk\tak&=0, \text{   $\Kt=1$}.
  \end{aligned}
\end{equation}
For the  ``diagonal" elements ($A<B,\/ K<L$):
\begin{multline}
  \tak\tbl-(\qabn)^{-1}\qkl(-1)^{\At\Lt+\Bt\Kt}\tbl\tak=\\
  (\qabn)^{-1}\frac{\l^{-\e-\h}-\l^{\e+\h}}{\l^{-\h}+\l^{\h}}
                                   (-1)^{(\At+\Bt)\Kt}\tbk\tal
\label{diag}
\end{multline}
For the  ``antidiagonal" elements ($A<B, K>L$):
\begin{multline}
  \tak\tbl-(\qabn)^{-1}\qkl(-1)^{\At\Lt+\Bt\Kt}\tbl\tak=\\
  (\qabn)^{-1}\frac{\l^{-\e+\h}-\l^{\e-\h}}{\l^{-\h}+\l^{\h}}
(-1)^{(\At+\Bt)\Kt}\tbk\tal \label{adiag}
\end{multline}
Notice that for $\e=\h$ the right-hand side of (\ref{adiag}) vanishes, and
for
$\e=-\h$ so do the right-hand side of (\ref{diag}). The factor in the
right-hand side
of (\ref{diag}) or (\ref{adiag}) respectively reduces to
$\l^{-\e}-\l^{\e}$ (in both cases).

Now, if one will fix the classical space \ $V,\quad {\dim V=2}$\ and
change only the quantization parameters
$\qab$, he will arrive to the example described in the introduction.

\newcommand{\Ai}{A_i}
\newcommand{\Aj}{A_j}
\newcommand{\Ak}{A_k}
\newcommand{\Ax}{A_x}
\newcommand{\Ay}{A_y}
\newcommand{\Az}{A_z}
\newcommand{\Ar}{A_r}
\newcommand{\As}{A_s}
\newcommand{\mijk}{{m_{ijk}}}
\newcommand{\mxar}{{m_{x\a r}}}
\newcommand{\mybs}{{m_{y\b s}}}
\newcommand{\Dxyz}{{\D_{xyz}}}
\newcommand{\Dkrs}{{\D_{krs}}}
\newcommand{\Dixy}{{\D_{ixy}}}
\newcommand{\Djab}{{\D_{j\a\b}}}

\section*{Appendix.}

The general definition of a {\it bialgebra $A=(\Ai)$ indexed
by the set}
$I$ is as follows. We should fix the relations $P\subset I^3,\/Q\subset
I^3$ and define the linear maps
\begin{align}
   \mijk:\Ai\otimes\Aj\lra\Ak  & \text{ for all } (i,j,k)\in P, \notag\\
   \Dxyz:\Ax\lra\Ay\otimes\Az  & \text{ for all } (x,y,z)\in Q, \notag
\end{align}
such that the diagram
$$
\begin{CD}
          \Ai\otimes\Aj        @>{\mijk}>> \Ak @>{\Dkrs}>>  \Ar \otimes \As \\
      @V{\Dixy\otimes\Djab}VV           @.          @AA{\mxar\otimes\mybs}A \\
     \Ax\otimes\Ay\otimes\Aa\otimes\Ab  @. \llra @.
                                           \Ax\otimes\Aa\otimes\Ay\otimes\Ab
\end{CD}
$$
is commutative when the maps in it make sense.

It is obvious that this definition
is self-dual: if ${(\Ai,I,P,Q,m,\D)}$ is a bialgebra, then
${(\Ai',I,Q,P,\D',m')}$ is also a bialgebra.
This definition contains in particular:
 graded and almost graded algebras, coalgebras and bialgebras (in the usual
sense), function algebras on categories,  additive categories,
function algebras on quantum categories,
algebras dual to them, etc.
Rather exotic
diagonals one can find, for example,
in the Odessky -- Feigin algebras~\cite{of}.

The possibility that need to be also considered
is that the indices
$i\in I$ may form not just a set but topological space or a smooth
manifold.
The natural farther step is to allow the indices
to be coordinates on a supermanifold or even a
 ``quantum manifold''.
In the context of this paper that means to allow to quantize the
quantization parameters themselves.

\end{document}